\documentclass[12pt,oneside]{amsart}
\usepackage{amsmath,amsfonts,latexsym,graphicx,amssymb,url}
\usepackage{amsmath,txfonts,pifont}
\usepackage{amsthm}
\newcommand{\trace}{\text{\rm trace\,}}
\newcommand{\diver}{\text{\rm div}}
\newcommand{\R}{{\mathbb R}} %%reals
\newcommand{\C}{{\mathbb C}}

\def \H {{\mathbb H}}
\def \e {\varepsilon}
\def \O {\Omega}
\def \l {\lambda}
\def \p {\partial}

\def \s {\sigma}

\def \R {{\mathbb R}}

\def \C {{\mathbb C}}

\def \s {{\sigma}}

\def \O {{\Omega}}

\def \p {{\partial}}

\def \l {{\lambda}}
\def  \e {{\varepsilon}}

\def \O {\Omega}

\def \a {\alpha}

\newtheorem{theorem}{Theorem}[section]
\newtheorem{corollary}[theorem]{Corollary}
\newtheorem{lemma}[theorem]{Lemma}
\newtheorem{proposition}[theorem]{Proposition}
\newtheorem{definition}[theorem]{Definition}
\newtheorem{remark}[theorem]{Remark}

\title[On the second order derivatives
of convex functions on the Heisenberg group]{On the second order derivatives
of convex functions on the Heisenberg group}
\author{
Cristian E. Guti\'errez and
Annamaria Montanari
 }
\thanks{\today\\
The first author was partially
supported by NSF grant DMS--0300004, and thanks the University of Bologna and the
INDAM for the kind hospitality and their support.
The second author was partially supported
by the University of Bologna, funds for selected research topics.}
\address{Department of Mathematics\\
    Temple University\\
    Philadelphia,
     PA 19122}

\email{gutierrez@math.temple.edu}

\address{Dipartimento di Matematica\\
    Universit\`a di Bologna\\
    Piazza Porta San Donato 5\\
   Bologna, 40127,  Italy}

\email{montanar@dm.unibo.it}
\begin{document}
\maketitle

%\printindex
%\centerline{\sc Preliminary Draft}

\setcounter{equation}{0}
\section{Introduction}
     \index{Introduction|(}

    % The above command produces an indexentry that is sorted
    % as "a--Raum", but has the appearance of "$\alpha--Raum$
    % in the makeindex output that is read by \printindex, i.e.
    % expands to
    %  "\item $\alpha--Raum$, <some "fill" cmd>  <Page number/range>
    % The "|(" opens a page-range that is closed by a "|)" on the

% last page of the range.

A classical result of Aleksandrov asserts that convex functions in
$\R^n$ are twice differentiable a.e., and a first step to prove it
is to show that these functions have second order distributional
derivatives which are measures, see \cite[pp. 239-245]{EG}.
%Such result is of great
%importance in the theory of Monge--Amp\`ere
%equations, and one way to prove it is to show that convex functions
%have second order distributional derivatives which are measures,
%see \cite[pp. 239-245]{EG}.
On the Heisenberg group, and more generally in Carnot groups,
several notions of convexity have been introduced and compared in
\cite{Danielli-Garofalo-Nhieu:notionsofconvexity} and
\cite{Lu-Manfredi-Stroffolini:notionsofconvexityinheisenberg}, and
Ambrosio and Magnani \cite[p. 3]{AM} ask the natural question if a
similar result holds in this setting.
%A reason for this question
%is that those estimates would be useful in the study of solutions
%for nondivergence equations of the form $a_{ij}X_i\,X_j$ where
%$a_{ij}$ is a uniformly elliptic measurable matrix and $X_i$ are
%the Heisenberg vector fields defined above.
Recently, these authors proved in \cite[Theorem 3.9]{AM} that
$BV^2_\H$ functions on Carnot groups, that is, functions whose
second order horizontal distributional derivatives are measures of
$H$-bounded variation, have second order horizontal derivatives
a.e., see Subsection \ref{subsec:ambrosiomagnani} below for
precise statements and definitions. On the other hand and also
recently, Lu, Manfredi and Stroffolini proved that if $u$ is an
$\mathcal H$--convex function in an open set of the Heisenberg
group $\H^1$ in the sense of the Definition \ref{def:Hconvex}
below, then the second order symmetric derivatives
$(X_iX_ju+X_jX_iu)/2,$ $i,j=1,2,$ are Radon measures \cite[Theorem
4.2]{Lu-Manfredi-Stroffolini:notionsofconvexityinheisenberg},
where $X_j$ are the Heisenberg vector fields defined by
\eqref{eq:heisenbergvectorfields}. Their proof is an adaptation of
the Euclidean one, it is based on the Riesz representation
theorem, and it can be carried out in the same way for $\H^n.$
However, to prove that $\mathcal H$--convex functions $u$ are
$BV^2_\H,$ one should show that the non symmetric derivatives
$X_iX_ju$ are Radon measures. Since the symmetry of the horizontal
derivatives is essential in the proof of \cite[Theorem
4.2]{Lu-Manfredi-Stroffolini:notionsofconvexityinheisenberg}, this
prevents these authors to answer the question of whether or not
the class of $\mathcal H$--convex functions is contained in
$BV^2_\H.$

The purpose in this paper is to establish the existence a.e. of
second order horizontal derivatives for the class of $\mathcal
H$-convex functions in the sense of Definition \ref{def:Hconvex}.
%on the Heisenberg group $\H^n$.
We will actually prove the stronger result that every $\mathcal
H$-convex function belongs to the class $BV^2_\H$ answering the
question posed by Ambrosio and Magnani in the setting of the
Heisenberg group. In order to do this we use the technique from
our work \cite{GM} which we shall briefly explain. Indeed,
following an approach recently used by Trudinger and Wang to study
Hessian equations \cite{Trudinger-Wang:hessianmeasuresI}, we
proved in \cite{GM} integral estimates in $\H^1$ in terms of the
following Monge--Amp\`ere type operator: $\det \mathcal H(u) +
12\,(u_t)^2$, see Definition \ref{def:Hconvex}. We first
established, by means of integration by parts, a comparison
principle for smooth functions, and then extended this principle
to ``cones". Together with the geometry in $\H^1$, this leads to
an Aleksandrov type maximum principle \cite[Theorem 5.5]{GM}.
Moreover, in \cite[Proposition 6.2]{GM} we proved an the estimate
of the oscillation of $\mathcal H$--convex functions. This
estimate furnishes $L^2$ estimates of the Lie bracket
$[X_1,X_2]u=-4\p_tu$ of $\mathcal H$--convex functions on $\H^1$
and permits to fill the gap between the results in \cite[Theorem
4.2]{Lu-Manfredi-Stroffolini:notionsofconvexityinheisenberg} and
\cite[Theorem 3.9]{AM}, and to prove that
\[
X_iX_ju=
\frac{[X_i,X_j]u}{2}+\frac{(X_iX_j+X_jX_i)u}{2}, \quad i,j=1,2,
\]
%\marginpar[left]{Can this\\sentence be\\clarified?}
are
Radon measures.

Following the route just described in $\H^1$, in this paper we
introduce in $\H^n$ the operator $\s_2(\mathcal H(u)) +12 n
u_t^2$, where $\s_2$ is the second elementary symmetric function
of the eigenvalues of the matrix $\mathcal H(u)$, we define the
notion of $\s_2(\mathcal H)$--convex function related to this
operator, and as a main tool we establish a comparison principle
for $\s_2(\mathcal H)$--convex functions, see Definition
\ref{def:sigmaHconvex} and Theorem \ref{thm:comparisonprinciple}.
In this frame, we next establish an oscillation estimate,
Proposition \ref{pro:osc}, which yields as a byproduct $L^2$
estimates of $\p_tu$ in $\H^n$ for a class of functions bigger
than the class of ${\mathcal H}$--convex functions. We apply these
estimates to obtain that the class of $\mathcal H$-convex
functions is contained in $BV^2_\H$, and as a corollary of
\cite[Theorem 3.9]{AM} it follows that $\mathcal H$-convex functions have
horizontal second derivatives a.e.

The paper is organized as follows. Section \ref{sec:preliminaries}
contains preliminaries about $\H^n$, $BV_\H$ functions, and the definitions of
$\mathcal H$--convexity and $\s_2(\mathcal H)$--convexity. In
Section \ref{sec:comparisonprinciple} we prove a comparison
principle for $C^2$ functions.  Section \ref{sec:Hmeasures}
contains the oscillation estimate and the construction of the
analogue Monge--Amp\`ere measures for $\s_2(\mathcal
H)$--convex functions. Finally, in Section \ref{sec:aleksandrovdiff} we prove
Aleksandrov's type differentiability theorem for $\mathcal H$--convex
functions in $\H^n.$

\medskip
 \noindent {\bf Acknowledgements.} We are greatly indebted to
  Bianca Stroffolini for some useful discussions.

At the
Workshop on
Second Order Subelliptic Equations and Applications, Cortona June 2003, we learnt
from Nicola Garofalo that in a joint paper with Federico
Tournier,
they extended to
higher dimensions the Aleksandrov type maximum principle proved by us in \cite[Theorem 5.5]{GM}
for ${\mathcal H}$--convex functions in $\H^1$.

\setcounter{equation}{0}
\section{Preliminaries, $\mathcal
H$--convexity and $\s_2(\mathcal
H)$--convexity}\label{sec:preliminaries}

%$\fint$

%\leftpointright

%\rightpointright

%\ding{230}

%\ding{43}

%\vskip 0.2in

%\largepencil

%\HandPencilLeft

Let $\xi=(x,y,t), \xi_0=(x_0,y_0,t_0)\in \R^n\times \R^n \times
\R$, and if $x=(x_1,\dots, x_n),$ $y=(y_1,\dots, y_n),$ then
$x\cdot y=\sum_{j=1}^nx_jy_j.$
The Lie algebra of $\H^n$ is
spanned by the left-invariant vector fields
\begin{equation}\label{eq:heisenbergvectorfields}
X_j=
\partial_{x_j} + 2y_j\,
\partial_t,\quad X_{n+j}=\partial_{y_j} - 2x_j\,
\partial_t \quad {\text{\rm
for}} \, j=1,\dots n.
\end{equation}
We have $[X_j,X_{n+j}]=X_jX_{n+j}-X_{n+j}X_j=- 4\partial_t$ for
every $j=1,\dots n,$ and $[X_j,X_{i}]=X_jX_{i}-X_{i}X_j=0$ for
every $i\neq n+j.$ If $\xi_0=(x_0,y_0,t_0)$,
then the non--commutative multiplication law in $\H^n$
is given by
\[
\xi_0\circ \xi=(x_0 + x,y_0 +y, t_0+t+2(x\cdot y_0-y \cdot x_0)),
\]
and we have $\xi^{-1}=-\xi$, $(\xi_0\circ
\xi)^{-1}=\xi^{-1}\circ  \xi_0^{-1}$. In $\H^n$ we define the
gauge function
\[
\rho(\xi)=\left((|x|^2+|y|^2)^2+t^2\right)^{1/4},
\]
and the distance
\begin{equation}\label{dist}
d(\xi,\xi_0)=\rho(\xi_0^{-1}\circ \xi).
\end{equation}
The group $\H^n$ has a family of dilations that are the group
homomorphisms, given by
\[
\delta_{\lambda}(\xi)=(\l x, \l y, \l^2 t)
\]
for  $\l>0.$
Then
\[
d(\delta_\l\xi,\delta_\l \xi_0)=\l \,d(\xi,\xi_0).
\]
For more details about $\H^n$ see \cite[Chapters XII and
XIII]{S:book}.

\subsection{$BV_\H$ functions}\label{subsec:ambrosiomagnani}
For convenience of the reader, we collect here some definitions and
a result from Ambrosio and Magnani \cite{AM} particularized to the Heisenberg group
that will be used in the proof of
Theorem \ref{thm:differ}.

We identify the vector field $X_j$ with the vector
$(e_j,\overrightarrow{0}, 2y_j)$ in $\R^{2n+1}$ for $j=1,\cdots
,n$, and with the vector $(\overrightarrow{0},e_j, -2x_j)$ for
$j=n+1,\cdots ,2n$. Here $e_j$ is the $j$th-coordinate vector in
$\R^n$ and $\overrightarrow{0}$ is the zero vector in $\R^n$.
Given $\xi=(x,y,t)\in \R^{2n+1}$, with this identification we let
$\{X_j(\xi)\}_{j=1}^{2n}$ be the vectors with origin at $\xi$ and
set $H_\xi=\text{span}\{X_j(\xi)\}$. The set $H_\xi$ is a
hyperplane in $\R^{2n+1}$. Given $\Omega\subset \R^{2n+1}$ we set
$H\Omega=\cup_{\xi\in \Omega}H_\xi$. Consider $\mathcal
T_{c,1}(H\Omega)$ the class functions $\phi:\Omega\to \R^{2n+1},$
$\phi=\sum_{j=1}^{2n}\phi_jX_j$ that are smooth and with compact
support contained in $\Omega$ and denote by $\|\phi
\|=\sup_{\xi\in \Omega} \sum_{j=1}^{2n}|\phi_j(\xi)|$.
%\marginpar[left]{I think \\that this \\is correct.\\ Do you
%agree\\
%with me?}
%\marginpar[left]{Is this\\correct?}\footnote{the class functions
%$\phi=\sum_{j=1}^{2n}\phi_jX_j$ that are smooth and with compact
%support contained in $\Omega$}
\begin{definition}
We say that the function $u\in L^1(\Omega)$ is of $H$-bounded variation if
\[
\sup\left\{ \int_\Omega u \,\diver_X \phi \,dx : \phi\in \mathcal T_{c,1}(H\Omega),\,
\|\phi \|\leq 1\right\}<\infty,
\]
%\marginpar[left]{What\\norm\\is $\|\phi \|$?} \footnote{I think
%that $\|\phi \|=\sup_{\xi\in \Omega} \sum_{j}|\phi_j(\xi)|$ or an
%equivalent norm}
where $\diver_X \phi =\sum_{i=1}^{2n}X_i\phi_i$.
The class of these functions is denoted by $BV_\H(\Omega)$.
\end{definition}

\begin{definition}
Let $k\geq 2$. The function $u:\Omega\to \R$ has $H$-bounded $k$
variation if the distributional derivatives $X_ju$, $j=1,\cdots,
2n$ are representable by functions of $H$-bounded $k-1$ variation.
If $k=1$, then $u$ has $H$-bounded $1$ variation if $u$ is of
$H$-bounded variation. The class of functions with $H$-bounded $k$
variation is denoted by $BV_\H^k(\Omega)$.
\end{definition}

\begin{theorem}[Ambrosio and Magnani \cite{AM}, Theorem 3.9]\label{thm:ambrosiomagnani}
If $u\in BV_\H^2(\Omega)$, then for a.e. $\xi_0$ in $\Omega$ there
exists a polynomial $P_{[\xi_0]}(\xi)$ with homogeneous degree
$\leq 2$ such that
\[
\lim_{r\to 0^+} \dfrac{1}{r^2}\,\fint_{U_{\xi_0,r}}
|u(\xi)-P_{[\xi_0]}(\xi)|\,d\xi =0,
\]
where $U_{\xi_0,r}$ is the ball centered at $\xi_0$ with radius
$r$ in the metric generated by the vector fields $X_j$, and
\[
P_{[\xi_0]}(\xi)=
P_{[\xi_0]}\left(\exp\left(\sum_{j=1}^{2n}\eta_jX_j+\eta_{2n+1}[X_1,X_2]\right)(\xi_0)\right)
=\sum_{|\a|\leq 2}c_\a \eta^\a,
\]
with $\a=(\a_1, \dots, \a_{2n+1}),$ $c_\a \in \R,$
$\eta^\a=\eta_1^{\a_1}\cdots \eta_{2n+1}^{\a_{2n+1}}$ and
$|\a|=\sum_{j=1}^{2n} \a_j+ 2 \a_{2n+1}.$
%\marginpar[left]{Can
%we\\explain in\\plain english\\this degree?}
\footnote{We can explicitly compute
$\eta=(x-x_0,y-y_0,(t_0-t+2(x\cdot y_0-y \cdot x_0))/4)$ by
solving the ODE $\xi=
\exp\left(\sum_{j=1}^{2n}\eta_jX_j+\eta_{2n+1}[X_1,X_2]\right)(\xi_0).$
}
\end{theorem}

\subsection{$\mathcal
H$--convexity and $\s_2(\mathcal H)$--convexity}

For a $C^2$ function $u$, let $X^2u$ denote the non symmetric matrix
$[X_iX_ju].$ Given $c\in \C$ and $u\in C^2(\Omega)$, let
\[
\mathcal H_c(u)= X^2u + c u_t \left[\begin{matrix} 0_n & I_n\\
-I_n & O_n\\
\end{matrix}\right].
\]
\begin{definition}\label{def:Hconvex}
The function $u\in C^2(\Omega)$ is $\mathcal H$--convex in
$\Omega$ if
the $2n\times 2n$ symmetric matrix
\[
\mathcal H(u)=\mathcal
H_2(u)=\left[\frac{X_iX_ju+X_jX_iu}{2}\right]
\]
is positive semidefinite in $\Omega$.
\end{definition}
Notice that the matrix $\mathcal H_c(u)$ is symmetric if and only
if $c=2.$ Also, if $\langle \mathcal H_c(u) \xi,\xi \rangle\geq 0$
for all $\xi\in \R^{2n}$ and for some $c$, then this quadratic form
is nonnegative for all values of $c\in \R$.

\begin{definition}\label{def:sigmaHconvex}
The function $u\in C^2(\Omega)$ is $\s_2(\mathcal H)$--convex in
$\Omega$ if
\begin{enumerate}
    \item the trace of the symmetric matrix $\mathcal H(u)$ is non negative,
    \item the second elementary symmetric function in the eigenvalues of $\mathcal H(u)$
\[
\s_2(\mathcal H(u))= \sum_{i<j}\left\{X^2_i u X^2_j u-
\left(\frac{X_iX_ju+X_jX_iu}{2}\right)^2 \right\}
\]
is non negative.
\end{enumerate}
\end{definition}

 We extend the definition of $\s_2(\mathcal
H)$--convexity to continuous functions.

\begin{definition}\label{def:convexity}
The function $u\in C(\Omega)$ is $\s_2(\mathcal H)$--convex in
$\Omega$ if there exists a sequence $u_k\in C^2(\Omega)$ of
$\s_2({\mathcal H})$--convex functions in $\Omega$ such that
$u_k\to u$ uniformly on compact subsets of $\Omega$.
\end{definition}

\begin{remark}\rm
If $u$ is $\mathcal H$--convex, then it is
$\s_2(\mathcal H)$--convex. The two definitions are equivalent in
$\H^1.$
Moreover, from \cite[Theorem 5.11]{Danielli-Garofalo-Nhieu:notionsofconvexity}
 we
have that if $u$ is convex in the standard sense, then $u$ is
$\mathcal H$--convex. However, the gauge function $\rho(x,y,t)=
\left((|x|^2+|y|^2)^2+t^2\right)^{1/4}$ is $\mathcal H$--convex but is
not convex in the standard sense.
\end{remark}

\section{Comparison Principle}\label{sec:comparisonprinciple}

A crucial step in the proof of Aleksandrov's type theorem, Theorem \ref{thm:differ}, is the following comparison principle for $C^2$ and
$\s_2(\mathcal H)$--convex functions.

\begin{theorem}\label{thm:comparisonprinciple}
Let $u,\varv\in C^2(\bar \Omega)$ such that $u+\varv$ is
$\s_2(\mathcal H)$--convex in $\Omega$ satisfying $\varv=u$ on
$\partial \Omega$ and $\varv<u$ in $\Omega.$ Then
\[
\int_\Omega \left\{ \s_2( \mathcal H(u))+ 12\,n\,\, (\partial_t u)^2
\right\}\,dz \leq \int_\Omega \left\{ \s_2(\mathcal H(\varv))+
12\,n\,\, (\partial_t \varv)^2 \right\}\,dz,
\]
and
\[
\int_\Omega \trace {\mathcal H}(u) \,dz \leq \int_\Omega \trace
{\mathcal H}(\varv) \,dz.
\]
\end{theorem}

\begin{proof}
By arguing as in \cite{GM}, set
\[
S(u) = \s_2(\mathcal H(u))= \sum_{i<j}\left\{X^2_i u X^2_j u-
\left(\frac{X_iX_ju+X_jX_iu}{2}\right)^2 \right\}.
\]
We have, by putting $r_{ij}=\dfrac{X_iX_ju+X_jX_iu}{2},$
%\marginpar[left]{I changed \\
%$r_{ij}$ to be\\
%consistent with\\
%Appendix.}
\begin{equation}\label{partial}
\dfrac{\partial S(u)}{\partial r_{ii}}=\sum_{j\neq i}X^2_j u;
\quad \dfrac{\partial S(u)}{\partial r_{ij}}=-\left(
\dfrac{X_iX_j+X_jX_i}{2} \right)u,
\end{equation}
and it is a standard fact that if $u$ is $\s_2({\mathcal
H})$--convex, then the matrix $\dfrac{\partial S(u)}{\partial
r_{ij}}$ is non negative definite, see Section \ref{sec:appendix}
for a proof.
%\marginpar[left]{Shall we \\ cite \cite{CS} too? }
% (see for example \cite{CS}).
%\marginpar[left]{Where\\is exactly\\the proof\\of this?}
%\footnote{pp. 264+277--277. I also have a simpler proof of it, and
%if you want we can add it in the Appendix.}
Let $0\leq s \leq 1$
and $\varphi(s)=S(\varv +s w), w=u-v$. Then
\begin{align*}
&\int_\Omega \{S(u)-S(\varv)\}\,dz\\ &= \int_0^1 \int_\Omega
\varphi'(s)\,dzds\\ &= \int_0^1 \int_\Omega \left\{
\sum_{i,j=1}^{2n} \frac{\partial S}{\partial r_{ij}}(\varv+sw)\,
(X_iX_j)w\right\}\,dzds
\\ &=
\int_0^1 \int_\Omega \left\{ \sum_{i,j=1}^{2n} X_i\left(
\frac{\partial S}{\partial r_{ij}}(\varv+sw) \,X_jw\right)
-
X_i\left( \frac{\partial S}{\partial r_{ij}}(\varv+sw)\right) X_jw
\right\}\,dzds\\ &= A- B.
\end{align*}
Since $w=0$ on $\p \Omega,$ $w>0$ in $\Omega,$ then the normal to $\p \Omega$ is
$\nu_X=-\dfrac{Xw}{|Dw|}$. Integrating by parts $A$ we have
\begin{align*}
A &= \int_0^1 \int_\Omega  \sum_{i,j=1}^{2n} X_i\left(
\frac{\partial S}{\partial r_{ij}}(\varv+sw)\right) \,X_jw \,dzds
\\
&= \int_0^1\int_{\partial\Omega} \sum_{i,j=1}^{2n} \left(
\frac{\partial S}{\partial r_{ij}}(\varv+sw) \right)\,X_jw
\,\nu_{X_i} d\s(z)ds\\
&=- \int_0^1\int_{\partial\Omega}
\sum_{i,j=1}^{2n} \left( \frac{\partial S}{\partial
r_{ij}}(\varv+sw) \,X_jw\right) \frac{{X_i}w}{|Dw|} d\s(z)ds\\
&=- \frac{1}{2}\int_{\partial\Omega}
\sum_{i,j=1}^{2n} \left( \frac{\partial S}{\partial
r_{ij}}(u+\varv) \,X_jw\right) \frac{{X_i}w}{|Dw|} d\s(z)\leq 0.
\end{align*}

We now calculate $B$. Let us remark that for any fixed $j=1,\dots,
2n$ by \eqref{partial} we have
\begin{align*}
\sum_{i=1}^{2n}X_i\left(\frac{\partial S}{\partial r_{ij}}\omega
\right)&= X_j\left(\frac{\partial S}{\partial r_{jj}}\omega
\right)+\sum_{i\neq j}X_i\left(\frac{\partial S}{\partial
r_{ij}}\omega \right)\\
&= X_j\left(\sum_{k\neq j}X^2_k\omega
\right)-\sum_{i\neq j}X_i\left(\frac{X_iX_j\omega +X_jX_i\omega }{2} \right)\\
&= \sum_{i\neq j}\left(X_jX^2_i\omega
-X_i\left(\frac{X_iX_j\omega +X_jX_i\omega }{2}
\right)\right)\\
&= \sum_{i\neq
j}\left(\frac{[X_j,X_i]X_i\omega}{2}+\frac{[X_j,X_i]X_i\omega}{2}+\frac{X_i[X_j,X_i]\omega}{2}
\right)\\
&= {3}\sum_{i\neq
j}\left(\frac{X_i[X_j,X_i]\omega}{2}
\right)\\
&= \frac{3}{2}
\begin{cases}
X_{j + n}[X_j,X_{j + n}]\omega, &\text{ if $j\leq n$}\\
X_{j - n}[X_j,X_{j - n}]\omega, &\text{ if $j > n$},
\end{cases}
\end{align*}
where, in the last two equalities, we have used the remarkable
fact that $[X_i,[X_j,X_k]]=0$ for every $i,j,k=1,\dots,2n,$ and
$[X_j,X_i]\ne 0$ iff $i=j\pm n.$
Hence,
\begin{align*}
B &= \int_0^1\int_\Omega \sum_{i,j=1}^{2n} X_i\left( \frac{\partial
S}{\partial r_{ij}}(\varv+sw)\right) X_jw\,dzds\\
&=\frac{3}{2}
\int_0^1\int_\Omega \sum_{j=1}^n X_{j + n}[X_j,X_{j + n}](\varv+sw) X_jw\,dzds\\
&\qquad \qquad \qquad + \frac32
\int_0^1\int_\Omega \sum_{j=n+1}^{2n} X_{j - n}[X_j,X_{j - n}](\varv+sw) X_jw\,dzds
\\
&=\frac{3}{2}
\int_0^1\int_\Omega \sum_{j=1}^n X_{j + n}\left\{ [X_j,X_{j + n}](\varv+sw)
X_jw\right\} \,dzds\\
&\qquad \qquad \qquad -\frac{3}{2}
\int_0^1\int_\Omega \sum_{j=1}^n [X_j,X_{j + n}](\varv+sw)
X_{j + n}X_jw\,dzds\\
&\quad +\frac{3}{2}
\int_0^1\int_\Omega \sum_{j=n+1}^{2n} X_{j - n}\left\{ [X_j,X_{j - n}](\varv+sw)
X_jw\right\} \,dzds\\
&\qquad \qquad \qquad -\frac{3}{2}
\int_0^1\int_\Omega \sum_{j=n+1}^{2n} [X_j,X_{j - n}](\varv+sw)
X_{j - n}X_jw\,dzds\\
&=\frac{3}{2}
\int_0^1\int_\Omega \sum_{j=1}^n X_{j + n}\left\{ -4\partial_t(\varv+sw)
X_jw\right\} \,dzds\\
&\qquad \qquad \qquad -\frac{3}{2}
\int_0^1\int_\Omega \sum_{j=1}^n [X_j,X_{j + n}](\varv+sw)
X_{j + n}X_jw\,dzds\\
&\quad +\frac{3}{2}
\int_0^1\int_\Omega \sum_{j=n+1}^{2n} X_{j - n}\left\{ 4\partial_t(\varv+sw)
X_jw\right\} \,dzds\\
&\qquad \qquad \qquad -\frac{3}{2}
\int_0^1\int_\Omega \sum_{j=n+1}^{2n} [X_j,X_{j - n}](\varv+sw)
X_{j - n}X_jw\,dzds\\
&=\frac{3}{2}
\int_0^1\int_\Omega \sum_{j=1}^n X_{j + n}\left\{ -4\partial_t(\varv+sw)
X_jw\right\} \,dzds\\
&\qquad \qquad \qquad -\frac{3}{2}
\int_0^1\int_\Omega \sum_{j=1}^n [X_j,X_{j + n}](\varv+sw)
X_{j + n}X_jw\,dzds\\
&\quad +\frac{3}{2}
\int_0^1\int_\Omega \sum_{j=1}^n X_j \left\{ 4\partial_t(\varv+sw)
X_{n+j}w\right\} \,dzds\\
&\qquad \qquad \qquad -\frac{3}{2}
\int_0^1\int_\Omega \sum_{j=1}^n [X_{j+n},X_j](\varv+sw)
X_jX_{j+n}w\,dzds
\end{align*}
\begin{align*}
&=\frac{3}{2}
\int_0^1\int_{\partial \Omega} \sum_{j=1}^n  -4\partial_t(\varv+sw)
X_jw \, \nu_{X_{j+n}} \,d\sigma(z)ds\\
&\qquad \qquad \qquad -\frac{3}{2}
\int_0^1\int_\Omega \sum_{j=1}^n [X_j,X_{j + n}](\varv+sw)
X_{j + n}X_jw\,dzds\\
&\quad +\frac{3}{2}
\int_0^1\int_{\partial \Omega} \sum_{j=1}^n   4\partial_t(\varv+sw)
X_{n+j}w\,\nu_{X_j} \,d\sigma(z)ds\\
&\qquad \qquad \qquad -\frac{3}{2}
\int_0^1\int_\Omega \sum_{j=1}^n [X_{j+n},X_j](\varv+sw)
X_jX_{j+n}w\,dzds\\
&=\frac{3}{2}
\int_0^1\int_{\partial \Omega} \sum_{j=1}^n  -4\partial_t(\varv+sw)
X_jw \, \nu_{X_{j+n}} \,d\sigma(z)ds\\
&\qquad \qquad \qquad -\frac{3}{2}
\int_0^1\int_\Omega \sum_{j=1}^n [X_j,X_{j + n}](\varv+sw)
[X_{j + n},X_j]w\,dzds\\
&\quad +\frac{3}{2}
\int_0^1\int_{\partial \Omega} \sum_{j=1}^n   4\partial_t(\varv+sw)
X_{n+j}w\,\nu_{X_j} \,d\sigma(z)ds\\
&=-\frac{3}{2}
\int_0^1\int_\Omega \sum_{j=1}^{n} [X_j,X_{j+ n}](\varv+sw)
[X_{j+ n},X_j]w\,dzds\\
&=\frac{3n}{2}
\int_0^1\int_\Omega (4\p_t)(\varv+sw)
(4\p_t)w\,dzds={24n}
\int_0^1\int_\Omega (\p_t\varv+s\p_tw)
\p_tw\,dzds\\
&= 12\,n\,\int_\Omega \{ (\partial_t u)^2 -(\partial_t\varv)^2\}\,dz.
\end{align*}
This completes the proof of the first inequality of the theorem.
The proof of the second one is similar.

\end{proof}

\section{Oscillation estimate and $\s_2({\mathcal H})$--Measures}\label{sec:Hmeasures}

In this section we prove that if $u$ is $\s_2({\mathcal
H})$--convex, we can locally control the integral of $\s_2( \mathcal
H)(u)+12\,n\,(u_t)^2$ in terms of the oscillation of $u.$ This
estimate will be crucial for the $L^2$ estimate of $\p_t u.$

Let us start with a lemma on $\s_2({\mathcal H})$--convex
functions.
\begin{lemma}\label{lem:comp}
If $u_1,u_2\in C^2(\Omega)$ are $\s_2({\mathcal H})$--convex, and
$f$ is convex in $\R^2$ and nondecreasing in each variable, then
the composite function $w=f(u_1,u_2)$ is $\s_2({\mathcal
H})$--convex.
\end{lemma}
\begin{proof}
Assume first that $f\in C^2(\R^2).$ We have
\[
X_jw=\sum_{p=1}^2\frac{\p f}{\p u_p}X_ju_p,
\]
\[
X_iX_jw=\sum_{p=1}^2\left(\frac{\p f}{\p u_p}X_iX_ju_p
+\sum_{q=1}^2\frac{\p^2 f}{\p u_q \p u_p}X_iu_qX_ju_p\right),
\]
and for every $h=(h_1,h_2)\in \R^2$
\[
\begin{split}
\langle \mathcal H(w) h,h\rangle &= \sum_{i,j=1}^{2n}X_iX_jw
\,h_i\,h_j
\\
&=\sum_{p=1}^2\frac{\p f}{\p u_p}\langle \mathcal H(u_p)
h,h\rangle +\sum_{p,q=1}^2\frac{\p^2 f}{\p u_q \p
u_p}(\sum_{i=1}^{2n}X_iu_qh_i)(\sum_{j=1}^{2n}X_ju_ph_j).
\end{split}
\]
Since the trace and the second elementary symmetric function of the
eigenvalues of the matrix $\mathcal H(u_p)$ are non negative,
$\dfrac{\p f}{\p u_p}\geq 0$ for $p=1,2$, and the matrix
\[
\left(\frac{\p^2 f}{\p u_q \p u_p}\right)_{p,q=1,2}
\]
is non negative definite, it follows that $w$ is $\s_2(\mathcal H)$--convex.

If $f$ is only continuous, then given $h>0$ let
\[
f_h(x)=h^{-2}\int_{\R^2} \varphi\left(\frac{x-y}{h}\right)f(y)dy,
\]
where $\varphi\in C^\infty$ is nonnegative vanishing outside the
unit ball of $\R^2,$ and $\int \varphi =1.$ Since $f$ is convex,
then $f_h$ is convex and by the previous calculation
$w_h=f_h(u_1,u_2)$ is $\s_2({\mathcal H})$--convex. Since $w_h\to
w$ uniformly on compact sets as $h\rightarrow 0$, we get that $w$
is $\s_2({\mathcal H})$--convex.
\end{proof}

\begin{proposition}\label{pro:osc}
Let $u\in C^2(\Omega)$ be $\s_2({\mathcal H})$--convex. For any
compact domain $\O'\Subset \O$ there exists a positive constant
$C$ depending on $\Omega'$ and $\Omega$ and independent of $u$,
such that
\begin{equation}
\label{eq:osc} \int_{\O'}\{\s_2( \mathcal H(u))+12\,n\,
(u_t)^2\}\,dz\leq C ({\rm osc}_\O u)^2.
\end{equation}
\end{proposition}
\begin{proof}
Given $\xi_0\in \O$ let $B_R=B_R(\xi_0)$ be a $d$--ball of radius
$R$ and center at $\xi_0$ such that $B_R\subset \O.$ Let $B_{\s
R}$ be the concentric ball of radius $\s R,$ with $0<\s<1.$
Without loss of generality we can assume $\xi_0=0,$ because the
vector fields $X_j$ are left invariant with respect to the
group of translations. Let $M=\max_{B_R} u$, then $u-M\leq 0$ in
$B_R$. Given $\varepsilon>0$ we shall work with the function
$u-M-\varepsilon <-\varepsilon$. In other words, by subtracting a
constant, we may assume $u<-\e$ in $B_R,$ for each given positive
constant $\e$ which will tend to zero at the end of the proof.

Define
\[
m_0=\inf_{B_R} u,
\]
and
\[
\varv(\xi)=\frac{m_0}{(1-\s^4)R^4}(R^4-\|\xi\|^4).
\]
Obviously $\varv=0$ on $\p B_R$ and $\varv=m_0$ on $\partial B_{\s
R}.$ We claim that $\varv$ is $\s_2({\mathcal H})$--convex in
$B_R$ and $\varv\leq m_0$ in $B_{\s R}.$ Indeed, setting $r=\|\xi\|^4$,
$h(r)=\dfrac{m_0}{(1-\s^4)R^4}(R^4-r)$, and following the
calculations in the proof of \cite[Proposition 6.2]{GM} we get
\[
\s_2 (\mathcal
H(\varv))=c_n(|x|^2+|y|^2)^2\left(\frac{m_0}{(1-\s^4)R^4}\right)^2\geq
0,
\]
with $c_n$ a positive constant and
 \[
\trace (\mathcal
H(\varv))=-(8n+4)\,(|x|^2+|y|^2)\frac{m_0}{(1-\s^4)R^4}\geq 0,
 \]
because $m_0$ is negative. Hence $\varv$ is $\s_2(\mathcal
H)$--convex in $B_R.$ Since $\varv-m_0=0$ on $\partial B_{\s R}$,
it follows from \cite[Proposition 5.1]{GM} that $\varv\leq m_0$ in
$B_{\s R}$. In particular, $\varv\leq u$ in $B_{\s R}.$

Let $\rho\in C_0^\infty(\R^2)$, radial with support in the
Euclidean unit ball, $\int_{\R^2}\rho(x)\,dx=1$, and let
\begin{equation}\label{eq:functionfh}
f_h(x_1,x_2)=h^{-2}\,
\int_{\R^2}\rho((x-y)/h)\,\max\{y_1,y_2\}\,dy_1dy_2.
\end{equation}

Define
\[
w_h=f_h(u,\varv).
\]
%Since the function $f=\max$ is convex, we can apply to it Lemma
%\ref{lem:comp}.
From Lemma \ref{lem:comp} $w_h$ is $\s_2(\mathcal H)$--convex in
$B_R$. If $y\in B_{\sigma R}$ then $\varv(y)\leq u(y)$. If
$\varv(y)<u(y)$ then $f_h(u,\varv)(y)=u(y)$ for $h$ sufficiently
small; and if $\varv(y)=u(y)$, then $f_h(u,\varv)(y)= u(y)+\alpha
\,h$. Hence
\begin{align}\label{eq:uinsigmarcontrolled}
\int_{B_{\s R}} \{\s_2( \mathcal H(u))+12\,n\,(\p_tu)^2\}\,dz &=
\int_{B_{\s R}}
\{\s_2( \mathcal H(w_h))+12\,n\,((w_h)_t)^2\}\,dz\notag\\
&\leq \int_{B_R} \{\s_2( \mathcal H(w_h))+12\,n\,((w_h)_t)^2\}\,dz.
\end{align}
Now notice that $f_h(u,\varv)\geq \varv$ in $B_R$ for all $h$
sufficiently small. In addition, $u<0$ and $\varv=0$ on $\partial
B_R$ so $f_h(u,\varv)=0$ on $\partial B_R$. Then we can apply
Theorem \ref{thm:comparisonprinciple} to $w_h$ and $\varv$ to get
\[
\begin{split}
\int_{B_R} \{\s_2(\mathcal H(w_h))+12\,n\,(\p_tw_h)^2\}\,dz &\leq
\int_{B_R}
\{\s_2( \mathcal H(\varv))+12\,n\,(\varv_t)^2\}\,dz\\
&= \left( \frac{m_0}{(1-\s)
R^4}\right)^2\int_{B_R}(c_n(|x|^2+|y|^2)^2+48n\,t^2)\,dz\\
& = \left( \frac{m_0}{(1-\s)
}\right)^2R^{2n-2}\int_{B_1}(c_n(|x|^2+|y|^2)^2+48n\,t^2)\,dz.
\\
\end{split}
\]
Combining this inequality with \eqref{eq:uinsigmarcontrolled} we get
\[
\int_{B_{\s R}} \{\s_2( \mathcal H(u))+12\,n\,(\p_tu)^2\}\,dz
\leq
C\,({m_0})^2R^{2n-2} \leq C\, R^{2n-2}({\rm osc}_{B_R} u +
\varepsilon)^2,
\]
and then \eqref{eq:osc} follows letting $\e\rightarrow
0$ and covering $\O'$ with balls.
\end{proof}
\begin{corollary}\label{corol}
Let $u\in C^2(\Omega)$ be $\s_2({\mathcal H})$--convex. For any
compact domain $\O'\Subset \O$ there exists a positive constant
$C,$ independent of $u,$ such that
\begin{equation}
\label{eq:osc2} \int_{\O'}\s_2( \mathcal H(u))\,dz\leq C ({\rm
osc}_\O u)^2R^{2n-2},
\end{equation}
and
\begin{equation}
\label{eq:osc21} \int_{\O'}(\p_tu)^2\,dz\leq C ({\rm osc}_\O
u)^2R^{2n-2}.
\end{equation}
\end{corollary}
\begin{corollary}\label{corol2}
Let $u\in C^2(\Omega)$ be $\s_2({\mathcal H})$--convex. For any
compact domain $\O'\Subset  \O$ there exists a positive constant
$C,$ independent of $u,$ such that
\begin{equation}
\label{eq:osc3} \int_{\O'}\trace {\mathcal H}_2(u)\,dz\leq C
R^{2n} {\rm osc}_\O u.
\end{equation}
\end{corollary}

\subsection{Measure generated by a $\s_2({\mathcal H})$--convex function}
We shall prove that the notion $\int \s_2( \mathcal H(u))+12\,n\,
\,u_t^2$ can be extended for continuous and $\s_2({\mathcal
H})$--convex functions as a Borel measure. We call this measure
the $\s_2({\mathcal H})$--measure associated with $u$, and we
shall show that the map $u\in C(\O)\rightarrow \mu(u)$ is weakly
continuous on $C(\O).$
\begin{theorem}\label{measure}
Given a $\s_2({\mathcal H})$--convex function $u\in C(\O)$, there
exists a unique Borel measure $\mu(u)$ such that when $u\in
C^2(\Omega)$ we have
\begin{equation}\label{eq:meas}
\mu(u)(E)=\int_E \{\s_2( \mathcal H(u))+12\,n\,u_t^2\}\,dz
\end{equation}
for any Borel set $E\subset \O.$ Moreover, if $u_k\in C(\Omega)$
are $\s_2({\mathcal H})$--convex, and $u_k\to u$ on compact
subsets of $\O,$ then $\mu(u_k)$ converges weakly to $\mu(u),$
that is,
\begin{equation}\label{eq:weakly}
\int_\O f\,d\mu(u_k) \rightarrow \int_\O f\,d\mu(u),
\end{equation}
for any $f\in C(\O)$ with compact support in $\O$.
\end{theorem}
\begin{proof} Let $u\in C(\O)$ be $\s_2({\mathcal H})$--convex, and let
$\{u_k\}\subset C^2(\O)$ be a sequence of $\s_2({\mathcal
H})$--convex functions converging to $u$ uniformly on compacts of
$\O$. By Proposition \ref{pro:osc}
\[
\int_{\O'}\{\s_2( \mathcal H(u_k))+12\,n \,(\p_t u_k)^2\}\,dz
\]
are uniformly bounded, for every $\O'\Subset\O$, and hence a
subsequence of $(\s_2( \mathcal H(u_k))+12n (\p_t u_k)^2)$
converges weakly in the sense of measures to a Borel measure
$\mu(u)$ on $\O.$
Moreover, by the same argument used in the proof of \cite[Theorem 6.5]{GM}
the map $u\in C(\Omega)
\rightarrow \mu(u)\in M(\O),$ the space of finite Borel measures
on $\O$, is well defined.

To prove \eqref{eq:weakly}, we first claim that it holds when
$u_k\in C^2(\Omega)$. Indeed, let $u_{k_m}$ be an arbitrary
subsequence of $u_k$, so $u_{k_m}\to u$ locally uniformly as $m\to
\infty$. By definition of $\mu(u)$, there is a subsequence
$u_{k_{m_j}}$ such that $\mu\left( u_{k_{m_j}} \right)\to \mu(u)$
weakly as $j\to \infty$. Therefore, given $f\in C_0(\Omega)$, the
sequence $\int_\Omega f\,d\mu(u_k)$ and an arbitrary subsequence
$\int_\Omega f\,d\mu(u_{k_m})$, there exists a subsequence
$\int_\Omega f\,d\mu(u_{k_{m_j}})$ converging to $\int_\Omega
f\,d\mu(u)$ as $j\to \infty$ and \eqref{eq:weakly} follows. For
the general case, given $k$ take $u_j^k\in C^2(\Omega)$
such that $u_j^k\to u_k$ locally uniformly as $j\to \infty,$
and then argue as
in the proof of \cite[Theorem 6.5]{GM}.
\end{proof}
\begin{corollary}\label{cor:Radon}
If $u \in C(\bar \Omega)$ is $\s_2(\mathcal H)$--convex in
$\Omega$, then $\mu(u)$ is a Radon measure.
\end{corollary}
\begin{proof}
The measure $\mu(u)$ is Borel regular, see \cite[p. 4-5]{EG} for definitions.
Indeed, given $A\subset \R^{2n+1}$ there exists open sets $V_k$ such that
$A\subset V_k$ and $\mu(u)(V_k)\leq \mu(u)(A)+1/k$ for all $k.$ Thus, $\mu(u)(A)=\mu(u)(\cap_1^\infty V_k).$
Finally, the estimate \eqref{eq:osc} implies that $\mu(u)(K)<\infty$ for all compact $K.$
Hence, $\mu(u)$ is a Radon measure.
\end{proof}
\begin{corollary}\label{cor:generalcomparison}
If $u,\varv\in C(\bar \Omega)$ are $\s_2(\mathcal H)$--convex in
$\Omega$, $u=\varv$ on $\partial \Omega$ and $u\geq \varv$ in
$\Omega$, then $\mu(u)(\Omega)\leq \mu(\varv)(\Omega)$.
\end{corollary}

By arguing as in \cite[Theorem 6.7]{GM} we also get the following
comparison principle for $\s_2(\mathcal H)$--measures.
\begin{theorem}\label{thm:comparisonpcpleformeasures}
Let $\Omega\subset \R^3$ be an open bounded set. If $u,\varv\in
C(\bar \Omega)$ are $\s_2(\mathcal H)$--convex in $\Omega$, $u\leq
\varv$ on $\partial \Omega$ and $\mu(u)(E)\geq \mu(\varv)(E)$ for
each $E\subset \Omega$ Borel set, then $u\leq \varv$ in $\Omega$.
\end{theorem}

\section{Aleksandrov-type differentiability theorem for $\mathcal H$--convex functions}\label{sec:aleksandrovdiff}
As an application of our previous results we finally have the following main theorem.
\begin{theorem}\label{thm:differ}
If $u$ is $\mathcal H$--convex, then $u\in BV_\H^2$ and so the
distributional derivatives $X_iX_ju$ exist a.e. for every
$i,j=1,\dots,2n.$
\end{theorem}
\begin{proof}
If $u$ is $\mathcal H$--convex, then by \cite[Theorem
3.1]{Lu-Manfredi-Stroffolini:notionsofconvexityinheisenberg}
$u$ is locally Lipschitz continuous with respect to the distance
$d$ defined in \eqref{dist},
 and $X_iu$ exists a.e. for $i=1,\dots, 2n.$  Moreover, by \cite[Theorem
4.2]{Lu-Manfredi-Stroffolini:notionsofconvexityinheisenberg} there is a
Radon measure $d\nu^{ij}$ such that, in the sense of distributions
\[
\frac{X_iX_ju+X_jX_iu}{2}=d\nu^{ij}, \quad i,j=1,\dots,2n.
\]
On the other hand, since $u$ is continuous and $\s_2(\mathcal
H)$--convex, then by \eqref{eq:osc21} $\p_tu$ is in $L^2_{loc}$.
Let $K\Subset \Omega,$ $\phi=\sum_{j=1}^{2n}\phi_jX_j
 \in C^2(\Omega, \R^{2n+1}),$
%\marginpar[left]{
%$\phi$ is not $\sum \phi_j e_j,$\\
%it is  $\sum\phi_jX_j$}
 with compact support in $K,$ $\|\phi\|<1.$
Since
\[
X_iX_j=\frac{X_iX_j+X_jX_i+[X_i,X_j]}{2}
=\frac{X_iX_j+X_jX_i}{2}\pm 2\delta_{i,i\mp n}\p_t,
\]
then for any $i=1,\dots, 2n$
\begin{equation}\label{eq:meas1}
\begin{split}
\int_\Omega X_iu \,\textrm{div}_X \,(\phi )dz=&-\int_\Omega u \,X_i\,\textrm{div}_X \, (\phi) dz\\
=&-\sum_{j=1}^{2n}\int_\Omega u X_iX_j\phi_j dz\\
=&-\sum_{j=1}^{2n}\int_\Omega u \left(\frac{X_iX_j\phi_j+X_jX_i\phi_j}{2}\pm 2\delta_{i,i\mp n}\p_t\phi_j\right) dz\\
=&\sum_{j=1}^{2n}\int_\Omega \phi_j d\nu^{ij}\mp 2\sum_{j=1}^{2n}\delta_{j\mp n,j}\int_\Omega u\,\p_t\phi_j dz\\
\leq & \sum_{j=1}^{2n}\nu^{ij}(K)\mp 2\sum_{j=1}^{2n}\delta_{j\mp n,j}\int_\Omega u\,\p_t\phi_j dz.
\end{split}
\end{equation}
Now, let $u_\e$ be the horizontal mollification of the function $u$
as in the proof of \cite[Theorem 4.2]{Lu-Manfredi-Stroffolini:notionsofconvexityinheisenberg}.
Then $u_\e$ is ${\mathcal H}$--convex and
\[
\left|\int_\Omega u_\e\p_t\phi_j dz\right|=\left|\int_\Omega \p_tu_\e\phi_j dz\right|
\leq c\|\p_tu_\e\|_{L^2(K)}\leq C,
\]
where  $c,C$ are positive constants depending on the diameter of $K$ and on the oscillation of $u$ over $K,$
but independent of $\e.$ Letting $\e$ tend to zero, we get
\begin{equation}\label{eq:meas2}
\left|\int_\Omega u \p_t\phi_j dz\right|\leq C.
\end{equation}
Thus, by \eqref{eq:meas1} and \eqref{eq:meas2} we can conclude that
\[
\int_\Omega X_iu \,\textrm{div}_X \,(\phi )dz\leq
\sum_{j=1}^{2n}\nu^{ij}(K)+C< \infty.
\]
Hence, $u\in BV_\H^2$ and the result then follows from Theorem
\ref{thm:ambrosiomagnani}.
\end{proof}

\small
\section{Appendix}\label{sec:appendix}
Let $A$ be an $n \times n$ symmetric matrix with eigenvalues
$\lambda_1,\cdots ,\lambda_n$, and the second elementary symmetric
function
\[
\s_2(A)=s(\l)=\sum_{j<k}\l_j\l_k
\]
with $\l=(\l_1,\dots, \l_n).$ An easy calculation shows that
\[
\frac{\p s}{\p \l_j}(\l)=\sum_{k\ne j}\l_k
\]
and
\begin{equation}\label{eq1}
s(\l)=\frac{1}{2}\left\{\left(\sum_{j=1}^n\l_j\right)^2-\sum_{j=1}^n\l_j^2\right\}.
\end{equation}
\begin{lemma}\label{lm:appdxbasic}
If $\s_2(A)\geq 0$ and $\text{\rm trace}(A)\geq 0$, then $\dfrac{\p s}{\p
\l_j}(\l)\geq 0$ for every $j=1,\dots, n.$
\end{lemma}
\begin{proof}
Since
\[
\text{\rm trace}(A)=\frac{\p s}{\p \l_j}(\l)+\l_j\geq 0,
\]
then either $\l_j\geq 0$ or $\dfrac{\p s}{\p \l_j}(\l)\geq 0.$ If
$\l_j\geq 0,$ since $s(\l)\geq 0,$ then by \eqref{eq1}
\[
\sum_{k=1}^n\l_k\geq \left(\sum_{k=1}^n\l_k^2\right)^{1/2}\geq
\l_j,
\]
and we get
\[
\frac{\p s}{\p \l_j}(\l)=\sum_{k\ne
j}\l_k=\sum_{k=1}^n\l_k-\l_j\geq 0.
\]
\end{proof}
\begin{proposition}
If $\s_2(A)\geq 0$ and $\text{\rm trace}(A)\geq 0$, then
\[\sum_{i,j=1}^n\frac{\p \s_2}{\p a_{ij}}(A)x_ix_j\geq 0
\]
for every $x\in \R^n.$
\end{proposition}
\begin{proof}
Let $C$ be a non negative definite Hermitian matrix. We write
\[
\s_2(A+C)-\s_2(A)=s(\eta_1,\dots,\eta_n)-s(\l_1,\dots,\l_n)
\]
where $\eta_1,\dots,\eta_n$ are the eigenvalues of $A+C$. Since
$C\geq 0,$ then $\eta_j\geq \l_j,$
 for any $j\in \{1,\dots,n\}.$ Moreover, by Lemma \ref{lm:appdxbasic},
 $\delta=\delta(A)=\dfrac{1}{2}\min\left\{
 \dfrac{\p s}{\p \l_j}(\l_1,\dots,\l_n) :\, j=1,\dots,n
\right\}\geq 0. $ If $C$ is small enough, then
\[\begin{split}
\s_2(A+C)-\s_2(A)
&=\int_0^1 \frac{d}{d\tau}s(\l+\tau(\eta-\l))\,d\tau\\
&=\sum_{j=1}^n\int_0^1 \frac{\p s}{\p
\l_j}(\l+\tau(\eta-\l))\,d\tau\,(\eta_j-\l_j)\\
&\geq  \delta
\sum_{j=1}^n (\eta_j-\l_j)=\delta \left({\text{\rm trace}}(A+C)-{\text{\rm trace}}(A)\right)
\\
&= \delta \, {\text{ \rm trace}} ( C)\geq 0.
\end{split}
\]
Let us now apply this inequality to the matrix
\[
C=tx\cdot x^T=t(x_ix_j), \quad x \in \R^n,
\]
and $t>0$ small enough. We obtain
\begin{equation}\label{101}
\s_2(A+tx \cdot x^T)-\s_2(A)\geq \delta \,{\text{\rm trace} \,}(C)=
\delta t |x|^2.
\end{equation}
On the other hand
\[
\frac{d}{dt}\s_2(A+tx \cdot x^T)\mid_{t=0}=\sum_{i,j=1}^n\frac{\p
\s_2}{\p a_{ij}}(A)x_ix_j.
\]
Then, from \eqref{101} we get
\begin{equation}\label{102}
\sum_{i,j=1}^n\frac{\p \s_2}{\p a_{ij}}(A)x_ix_j \geq \delta
|x|^2\geq 0, \quad \forall x\in \R^n.
\end{equation}
\end{proof}


\normalsize

\begin{thebibliography}{111}
\bibliographystyle{alpha}
\bibitem{AM} {L. Ambrosio and V. Magnani,} {\it Weak diferentiability of BV functions on stratified groups.}
\url{http://cvgmt.sns.it/papers/ambmag02/}

%\bibitem{CS} {L. A. Caffarelli, L. Nirenberg, and J. Spruck,} {\it
%The Dirichlet problem for non-linear second order elliptic
%equations III: Functions of the eigenvalues of the Hessian.} Acta
%Math. {\bf 155} (1985), 261-301.

\bibitem{Danielli-Garofalo-Nhieu:notionsofconvexity} D. Danielli,
N. Garofalo, and D. M. Nhieu, {\it Notions of convexity in Carnot
groups.}
Comm. in Analysis and Geometry,
{\bf 11} 2 (2003)
 263--341.

\bibitem{EG} {L. C. Evans and R. Gariepy,} {\it Measure theory and fine properties of functions,}
Studies in Advanced Mathematics, CRC Press, Boca Raton 1992.

%\bibitem{GTo}{N. Garofalo, F. Tournier,} Work in preparation.





\bibitem{GM} C. E. Guti\'errez and A. Montanari, {\it Maximum and comparison principles on the Heisenberg
group}. ArXiv preprint 2003.


\bibitem{Lu-Manfredi-Stroffolini:notionsofconvexityinheisenberg}
Guozhen Lu, J. Manfredi, and B. Stroffolini, {\it Convex functions on the
Heisenberg group.} To appear in Calculus of Variations.

\bibitem{S:book} E. M. Stein, {\it Harmonic Analysis: Real
Variable methods, Orthogonality and Oscillatory Integrals.} Vol. 43
of the Princeton Math. Series. Princeton U. Press. Princeton, NJ,
1993.



\bibitem{Trudinger-Wang:hessianmeasuresI} N. S. Trudinger and Xu-Jia
Wang, {\it Hessian measures I.} Topol. Methods Nonlinear Anal. {\bf
10} (1997) 225-239.




\end{thebibliography}
\end{document}